\documentclass[dvips,ejs]{imsart}

\usepackage{amsmath}
\usepackage{amssymb}
\usepackage{amsfonts}
\usepackage{amsthm}
\usepackage{rotating}
\usepackage{mathtools}

\usepackage{array}
\usepackage{multirow}
\usepackage{multicol}

\usepackage{graphicx}

\usepackage{tabularx}
\usepackage{longtable}
\usepackage{lscape}
\usepackage{float}

 \usepackage{hyperref}

\doi{10.1214/154957804100000000}
\pubyear{0000}
\volume{0}
\firstpage{0}
\lastpage{0}

\startlocaldefs
\newtheorem{prop}{Proposition}
\newtheorem{cor}{Corollary}
\newtheorem{theor}{Theorem}
\newtheorem{lem}{Lemma}
\newtheorem{defi}{Definition}

\newenvironment{preuve}[2]{
\noindent\itshape{#1} \\  
\small\itshape{#2}}
{ 
\vspace{0.25cm}\hspace{\stretch{3}}
\rule{1ex}{1ex}
}

\newcounter{moncompt} 
\renewcommand{\themoncompt}{\arabic{moncompt}} 

\newenvironment{procedure}[1]{
\addtocounter{moncompt}{0}
\refstepcounter{moncompt} 

{} 
\bigskip 
\emph{\textbf{ Procedure \themoncompt}}
\vspace{3mm}
\itshape{ #1}
}

\newcounter{mescond} 
\renewcommand{\themescond}{C.\arabic{mescond}} 

\newenvironment{cond}[1]{ 
\addtocounter{mescond}{0}
\refstepcounter{mescond} 
{}
\vspace{0.5cm}
\begin{center}
 \framebox[1.05\width]{\textit{#1} }\quad (\themescond) 
\end{center}
\vspace{0.5cm}
}

\newcommand{\id}[1]{
\text{#1}
}

\newcounter{moncomptexp} 
\renewcommand{\themoncomptexp}{\arabic{moncomptexp}} 

\newenvironment{exple}[1]{
\addtocounter{moncomptexp}{0}
\refstepcounter{moncomptexp} 

{} 
\vspace{3mm} 
\textit{\textbf{ Example \themoncomptexp:}}{#1}
}

\newlength{\eqboxstorage}
\newcommand{\eqbox}[1]{
  \setlength{\eqboxstorage}{\fboxsep}
  \setlength{\fboxsep}{6pt}
  \boxed{#1}
  \setlength{\fboxsep}{\eqboxstorage}
}

\endlocaldefs

\begin{document}

\begin{frontmatter}

\title{Generalized Hoeffding-Sobol Decomposition for Dependent Variables - Application to Sensitivity Analysis}
\runtitle{Generalized ANOVA Decomposition for Dependent Variables}


\author{\fnms{Gaelle} \snm{Chastaing}\ead[label=e1]{gaelle.chastaing@imag.fr}}
\address{Universit\'e de Grenoble, LJK/MOISE BP 53, 38041 Grenoble Cedex, France\\ \printead{e1}}

\author{\fnms{Fabrice} \snm{Gamboa}\ead[label=e2]{fabrice.gamboa@univ-tlse.fr}}
\address{Universit\'e Paul Sabatier, IMT-EPS,  118, Route de Narbonne, 31062 Toulouse  Cedex 9, France\\ \printead{e2}}
\author{\fnms{Cl\'ementine} \snm{Prieur}\ead[label=e3]{clementine.prieur@imag.fr}}
\address{Universit\'e Joseph Fourier, LJK/MOISE BP 53, 38041 Grenoble Cedex, France\\ \printead{e3}}



\runauthor{G. Chastaing {\it{et al.}}}


\begin{abstract}
In this paper, we consider a regression model built on dependent variables. This regression modelizes an input output relationship. Under boundedness type assumptions on the joint distribution function of the input variables, we show that a generalized Hoeffding-Sobol decomposition is available. This leads to new indices measuring the sensitivity of the output with respect to the input variables. We also study and discuss the estimation of these new indices.
\end{abstract}




\tableofcontents

\end{frontmatter}

\section {Introduction}

Sensitivity analysis (SA) aims to identify the variables that most contribute to the variability into a non linear regression model. 
Global SA is a stochatic approach whose objective is to determine a global criterion based on the density of the joint probability distribution function of the output 
and the inputs of the regression model. The most usual quantification is the variance-based method, widely studied in SA literature. 
Hoeffding decomposition~\cite{hoef} (see also Owen \cite{owen2}) states that the variance of the output can be uniquely decomposed into summands of increasing dimensions under orthogonality constraints. Following this approach, Sobol~\cite{sobol} introduces variability measures, the so called \textit{Sobol sensitivity indices}. These indices quantify the contribution of each input on the system.\\

Different methods have been exploited to estimate Sobol indices. The Monte Carlo algorithm was proposed by Sobol~\cite{sobol2}, and has been later improved by the Quasi Monte Carlo technique, performed by Owen~\cite{owen}. 
FAST methods are also widely used to estimate Sobol indices. 
Introduced earlier by Cukier {\it{et al.}}~\cite{cukier2} \cite{cukier}, they are well known to reduce the computational cost of multidimensional integrals thanks to Fourier transformations. 
Later, Tarantola {\it{et al.}}~\cite{tarantule} adapted the Random Balance Designs (RBD) to FAST method for SA (see also recent advances on the subject by Tissot {\it{et al.}}~\cite{tissot}).

However, these indices are constructed on the hypothesis that input variables are independent, which seems unrealistic for many real life phenomena. 
In the literature, only a few methods and estimation procedures have been proposed to handle models with dependent inputs. Several authors have proposed sampling techniques to compute marginal contribution of inputs to the outcome variance (see the introduction in Mara and references therein \cite{mara}). 
As underlined in Mara {\it{et al.}}\cite{mara}, if inputs are not independent, the amount of the response variance due to a given factor may be influenced by its dependence to other inputs. Therefore, classical Sobol indices and FAST approaches for dependent variables are difficult to interpret (see, for example, Da Veiga's illustration~\cite{daveiga1} p.133). 
Xu and Gertner~\cite{xu} proposed to decompose the partial variance of an input into a correlated part and an uncorrelated one. Such an approach allows to exhibit inputs that have an impact on the output only through their strong correlation with other incomes. However, they only investigated linear models with linear dependences.

Later, Li {\it{et al.}}~\cite{li} extended this approach to more
general models, using the concept of High Dimensional Model Representation (HDMR~\cite{rshdmr}).  
HDMR is based on a hierarchy of component functions of increasing dimensions (truncation of Sobol decomposition in the case of independent variables). 
The component functions are then approximated by expansions in terms of some
suitable basis functions (e.g., polynomials, splines, ...). This meta-modeling approach allows the splitting of
the response variance into a correlative contribution and a structural one of a set of inputs. 
Mara {\it{et al.}}~\cite{mara} proposed to decorrelate the inputs with the Gram-Schmidt
procedure, and then to perform the ANOVA-HDMR of Li {\it{et al.}}~\cite{li} on
these new inputs. The obtained indices can be interpreted as fully, partially correlated
and independent contributions of the inputs to the output. Nevertheless,
this method does not provide a unique orthogonal set of inputs as it depends on
the order of the inputs in the original set. Thus, a large number of sets has to be
generated for the interpretation of resulting indices. 
As a different approach, Borgonovo {\it{et al.}} ~\cite{bor,bor2} initiated the construction of a
new generalized moment free sensitivity index. Based on geometrical consideration,
these indices measure the shift area between the outcome density and this
same density conditionally to a parameter. Thanks to the properties of these
new indices, a methodology is given to obtain them analytically through test
cases.

Notice that none of these works has given an exact and unambiguous definition of the functional ANOVA for correlated inputs as the one provided by Hoeffding-Sobol decomposition when inputs are independent. Consequently, the exact form of the model has neither been exploited to provide a general variance-based sensitivity measures in the dependent frame. 

In a pionnering work, Hooker~\cite{hooker}, inspired by Stone~\cite{stone}, shed new lights on hierarchically orthogonal function decomposition. 
We revisit and extend here the work of Hooker. We obtain hierarchical functional decomposition under a general assumption on the inputs distribution. 
Furthermore, we also show the uniqueness of the decomposition leading to the definition of new sensitivity indices. 
Under suitable conditions on the joint distribution function of the input variables, we give a hierarchically orthogonal functional decomposition (HOFD) of the model.
The summands of this  decomposition are functions depending only on a subset of input variables and are hierarchically uncorrelated. This means that two of these components are orthogonal whenever all the variables involved in one of the summands also appear in the other. This decomposition leads to the construction of generalized sensitivity indices well tailored to perform global SA  when the input variables are dependent.   
In the case of independent inputs, this decomposition is nothing more than the Hoeffding one. Furthermore, our generalized  sensitivity indices are in this case the classical Sobol ones. In the general case,  the computation of the summands of the HOFD involves a minimization problem under constraints (see Proposition \ref{dugenou}). 
A statistical procedure to approach the solution of this counstrainted optimization problem will be investigated in a next paper.
Here, we will focus on the particular case where the inputs are independent pairs of dependent variables (IPDV).
Firstly, in the simplest case of a single pair of dependent variables, the HOFD may be performed by solving a functional linear system
of equations involving suitable projection operators (see Procedure \ref{proc1}). In the
more general IPDV case, the HOFD is then obtained in two steps (see Procedure
\ref{proc2}). The first step is a classical Hoeffding-Sobol decomposition of the output on
the input pairs, as developped in Jacques {\it{et al.}}~\cite{jac}. The second step is the
HOFDs of all the pairs. In practical situations, the non parametric regression
function of the model is generally not exactly known. In this case, one can
only have at hand some realizations of the model and have to estimate, with
this information, the HOFD. Here, we study this statistical problem in the
IPDV case. We build estimators of the generalized sensitivity indices and study
numerically their properties.  One of the main conclusion is that the
generalized indices have a total normalized sum. This is not true for classical
Sobol indices in the frame of dependent variables.
$\;$
\\
The paper is organized as follows.\\
In Section \ref{sect2}, we give and discuss general results on the HOFD. The main result
is Theorem \ref{theor1}. We show here that a HOFD is available under a boundedness type
assumption (\ref{c2}) on the density of the joint distribution function of the inputs. Further,
we introduce the generalized indices.
In Section \ref{sect3}, we give examples of multivariate distributions to which Theorem \ref{theor1} applies. We also state a sufficient condition for (\ref{c2}) and necessary and sufficient conditions in the IDPV case.
Section \ref{sect4} is devoted to the estimation procedures of the components of the HOFD and of the new sensitivity indices. 
Section \ref{sect5} presents numerical applications.
Through three toy functions, we estimate generalized indices and compare their
performances with the analytical values. In Section \ref{sect6}, we give conclusions and
discuss future work. Technical proofs and further details are postponed to the
Appendix.
%

\section{Generalized Hoeffding decomposition-Application to SA}\label{sect2}

To begin with, let introduce some notation. We briefly recall the usual functional ANOVA decomposition, and Sobol indices. We then state a generalization of this decomposition, allowing to deal with correlated inputs.

\subsection{Notation and first assumptions}
We denote by $\subset$ the strict inclusion, that is $A \subset B \Rightarrow A \cap B \neq B$, whereas we use $\subseteq$ when equality is possible.\\

Let $(\Omega,\mathcal{A},P)$ be a probability space and let $Y$ be the output of a deterministic model $\eta$.
Suppose that $\eta$ is a function of a random vector $\mathbf X=(X_1,\cdots,X_p) \in \mathbb R^p$, $p\geq 1$ and that $P_X$ is the pushforward measure of $P$ by $\mathbf X$, 

 $$
 \begin{array}{cc}
 Y :&
\begin{array}{ccccc}
  (\Omega,\mathcal{A},P)& \rightarrow &(\mathbb{R}^p,\mathcal{B}(\mathbb{R}^p),P_X) &\rightarrow & (\mathbb{R},\mathcal{B}(\mathbb{R})) \\
\omega & \mapsto & \mathbf X(\omega) & \mapsto & \eta(\mathbf X(\omega))
\end{array}

\end{array}
$$

 Let $\nu$ be a $\sigma$\textendash finite measure on $(\mathbb{R}^p,\mathcal{B}(\mathbb{R}^p))$. Assume that $P_X << \nu$ and let $p_X$ be the density of $P_X$ with respect to $\nu$, that is $p_X=\dfrac{dP_X}{d\nu}$.\\
Also, assume that $\eta \in L^2_{\mathbb R}(\mathbb R^p,\mathcal{B}(\mathbb{R}^p),P_X)$. The associated innner product of this Hilbert space is:
\[
 \langle h_1,h_2\rangle=\int h_1(\mathbf x)h_2(\mathbf x) p_X d\nu(\mathbf x) = \mathbb E(h_1(\mathbf X)h_2(\mathbf X)) 
\]

Here $\mathbb E(\cdot)$  denotes the expectation. The corresponding norm will be classically denoted by $\|\cdot\|$.\\ 
Further, $V(\cdot)=\mathbb E[(\cdot-\mathbb E(\cdot))^2]$ denotes the variance, 
and $\mbox{Cov}(\cdot,\ast)=\mathbb E[(\cdot-\mathbb E(\cdot))(\ast-\mathbb E(\ast))]$ the covariance.\\

Let $\mathcal P_p : =\{1,\cdots,p\}$ and $S$ be the collection of all subsets of $\mathcal P_p$.\\
Define $S^- := S\setminus \mathcal P_p$ as the collection of all subsets of $\mathcal P_p$ except $\mathcal P_p$ itself.\\

Further, let $X_u :=(X_l)_{l \in u}$, $u \in S\setminus \{\emptyset\}$. 
We introduce the subspaces of $L^2_{\mathbb R}(\mathbb R^p,\mathcal{B}(\mathbb{R}^p),P_X)$ $(H_u)_{u \in S}$, $(H_u^0)_{u \in S}$ and $H^0$ . $H_u$ is the set of all measurable and square integrable functions depending only on $X_u$. $ H_{\emptyset}$ is the set of constants and is identical to $(H_{\emptyset}^0)_{u \in S}$. $H_u^0$, $u \in S\setminus \emptyset$, and $H^0$ are defined as follows:

\[
 H_u^0=\left\{ h_u(X_u)\in H_u, ~\langle h_u,h_v\rangle=0, \forall~v \subset u, \forall~h_v \in H_v^0 \right\}
\]

\[
 H^0=\left\{ h(X)=\sum_{u \in S}h_u(X_u), h_u \in H_u^0\right\}
\]

At this stage, we do not make assumptions on the support of $\mathbf X$. For $u \in S\setminus \emptyset$, the support of $X_u$ is denoted by $\mathcal X_u$.

\subsection{Sobol sensitivity indices}\label{sect22}

In this section, we recall the classical Hoeffding-Sobol decomposition, and the Sobol sensitivity indices if the inputs are independent, that is when $P_X=P_{X_1} \otimes \cdots \otimes P_{X_p}$. \\
The usual presentation is done when $\mathbf X \sim \mathcal U([0,1]^p)$ \cite{sobol}, but the Hoeffding decomposition remains true in general case~\cite{vandervaart}.\\

Let $\mathbf x=(x_1,\cdots,x_p) \in \mathbb R^p$ and assume that $\eta \in \mathbb{L}^2(\mathbb R^p,P_X)$.
The decomposition consists in writting $\eta(\mathbf{x})=\eta(x_1,\cdots,x_p)$ as the sum of increasing dimension functions:

\begin{eqnarray}\label{dec}
 \eta(\mathbf{x}) &= &\eta_0+\sum_{i=1}^p \eta_i(x_i)+ \sum_{1\leq i <j \leq p}\eta_{i,j}(x_i,x_j)+ \cdots+ \eta_{1,\cdots,p}(\mathbf{x}) \nonumber\\
&=& \sum_{u \subseteq \{1 \cdots p\}} \eta_u (x_u)
\end{eqnarray}
\

 The expansion (\ref{dec}) exists and is unique under one of the hypothesis: 
$$
\left\{
\begin{array}{ccc}\label{cond}
\textrm{i)} & \int \eta_u(x_u) dP_{X_i} = 0 &\quad \forall ~i \in u, \forall ~u \subseteq \{ 1\cdots p \}\\
&\textrm{ or } &\\
 \textrm{ii)} & \int \eta_u(x_u) \eta_v(x_v) dP_X =0 &\quad \forall~ u,v \subseteq \{1 \cdots p \}, \quad u \neq v
\end{array}
\right.
$$

Equation (\ref{dec}) tells us that the model function $Y=\eta(\mathbf X)$ can be expanded in a functional ANOVA. 
The independence of the inputs and the orthogonality properties ensure the global variance decomposition of the output as $ V(Y)=\sum_{u \in S}V(\eta_u(X_u))$.

Moreover, by integration, each term $\eta_u$ has an explicit expression, given by:

\begin{equation}\label{sd}
 \eta_0=\mathbb E(X), \quad \eta_i=\mathbb E(Y/X_i)-\mathbb E(Y), ~ i=1,\cdots,p,\quad \eta_u=\mathbb E(Y/X_u)-\sum_{v \subset u} \eta_v,~ |u|\geq 2
\end{equation}

Hence, the contribution of a group of variables $X_u$ in the model can be quantified in the fluctuations of $Y$. The Sobol indices expressions are defined by:

\begin{equation}\label{eq1}
 S_u=\dfrac{V(\eta_u)}{V(Y)}=\dfrac{V[\mathbb E(Y/X_u)]-\sum_{v \subset u} V[\mathbb E(Y/X_v)]}{V(Y)}, \quad u \subseteq \mathcal P_p
\end{equation}
Furthermore,
\[
 \sum_{u \in S} S_u=1
\]

\bigskip

However, the main assumption is that the input parameters are independent. This is unrealistic in many cases. 
The use of expressions previously set up is not excluded in case of inputs' dependence, 
but they could lead to an unobvious and sometimes a wrong interpretation. 
Also, technics exploited to estimate them may mislead final results because most of them are built on the hypothesis of independence. 
For these reasons, the objective of the upcoming work is to show that the construction of sensitivity indices under dependence condition can be done into a mathematical frame. \\

In the next section, we propose a generalization of the Hoeffding decomposition under suitable conditions on the joint distribution function of the inputs. 
This decomposition consists of summands of increasing dimension, like in Hoeffding one. 
But this time, the components are hierarchically orthogonal instead of being mutually orthogonal. 
The hierarchical orthogonality will be mathematically defined further.
Thus, the global variance of the output could be decomposed as a sum of covariance terms depending on the summands of the HOFD. 
It leads to the construction of generalized sensitivity indices summed to $1$ to perform well tailored SA in case of dependence.

\subsection{Generalized decomposition for dependent inputs}
We no more assume that $P_X$ is a product measure. Nevertheless, we assume: 

\cond{
$
 \begin{array}{lc}
&P_X << \nu \\
\textrm{where}&\\
&\nu(dx)=\nu_1(dx_1) \otimes \cdots \otimes \nu_p(dx_p)
\end{array}
$
}\label{c1}

Our main assumption is :
\cond{
$
\begin{array}{lll}
\exists~ 0<M\leq 1 , ~\forall~u \subseteq \mathcal P_p,&  p_X\geq M \cdot p_{X_u}p_{X_{u^c}} \quad \nu\textrm{-a.e.} 
\end{array}
$
 }\label{c2}

 where $u^c$ denotes the complement set of $u$ in $\mathcal P_p$. $ p_{X_u}$ and  $p_{X_{u^c}}$ are respectively the marginal densities of $X_u$ and $X_{u^c}$.\\

The section is organized as follows: a preliminary lemma gives the main result to show that $H^0$ is a complete space. 
Then, this ensures the existence and the uniqueness of the projection of $\eta$ onto $H^0$.
The generalized decomposition of $\eta$ is finally obtained by adding a residual term orthogonal to every summand, as suggested in \cite{hooker}.
The first part of the reasoning is mostly inspired by Stone's work~\cite{stone}, except that our assumptions are more general.
Indeed, we have a relaxed condition on the inputs distribution function.
Moreover, the support $\mathcal X$ of $\mathbf X$ is general.\\

To begin with, let us state some definitions.
In the usual ANOVA context, a model is said to be hierarchical if for every term involving some inputs, all lower-order terms
involving a subset of these inputs also appear in the model. Correspondingly, a hierarchical collection $T$ of subsets of $\mathcal P_p$ is defined as follows:

\begin{defi}
A collection $T \subset S$  is hierarchical if for $u \in T$ and $v$ a subset of $u$, one has $v \in T$.
\end{defi}

The next Lemma is a generalization of the Lemma 3.1 of \cite{stone}. As already mentioned, it will be the key to show the hierarchical decomposition. 

\begin{lem}\label{lemstone}
Let $T \subset S$ be hierarchical. Suppose that (\ref{c1}) and (\ref{c2}) hold. Set $\delta=1-\sqrt{1-M} \in ]0,1]$. Then, for any $h_u \in H_u^0$, $u \in T$, we have:

\begin{equation}\label{lem}
 \mathbb E[(\sum_{u \in T }h_u(\mathbf X))^2] \geq \delta^{\#(T)-1} \sum_{u \in T} \mathbb E[h_u^2(\mathbf X)]
\end{equation}

\end{lem}

The proof of Lemma  \ref{lemstone} is postponed to the Appendix. Our main theorem follows:

\begin{theor}\label{theor1}
Let $\eta$ be any function in $L^2_{\mathbb R}(\mathbb R^p,\mathcal{B}(\mathbb{R}^p),P_X)$. Then, under (\ref{c1}) and (\ref{c2}), there exist functions $\eta_0,\eta_1,\cdots,\eta_{\mathcal P_p} \in H_{\emptyset}\times H_1^0 \times \cdots H_{\mathcal P_p}^0$ such that the following equality holds :
\begin{eqnarray}\label{equ2}
  \eta(X_1,\cdots,X_p)&=&\sum_i \eta_i(X_i)+\sum_{i,j}\eta_{ij}(X_i,X_j)+\cdots+\eta_{\mathcal P_p}(X_1,\cdots,X_p)\nonumber\\
&=&\sum_{u \in S} \eta_u(X_u)
\end{eqnarray}

Moreover, this decomposition is unique.
\end{theor}

The proof is given in the Appendix.\\
Notice that, in case where the input variables $X_1,\cdots,X_p$ are independent,  $\delta=1$ and Inequality (\ref{lem}) of Lemma \ref{lemstone} is an equality. Indeed, in this case, this equality is directly obtained by orthogonality of the summands.\\

The variational counterpart of Theorem \ref{theor1} is a minimization problem under conditional constraints. 

\begin{prop}\label{dugenou}
Suppose that (\ref{c1}) and (\ref{c2}) hold. Let $(\mathcal P)$ be the minimization problem under constraints:

$$
(\mathcal P)\left\{
\begin{array}{l}
 \displaystyle{\min_{(\tilde\eta_u)_{u \in S}} ~\mathbb E [(Y-\sum_{u \in S}\tilde\eta_u(X_u))^2]}\\
 \displaystyle{\mathbb E(\tilde\eta_u(X_u)/ X_{u\setminus i})=0, ~\forall~i \in u, ~\forall~u \in S\setminus \emptyset}
\end{array}
\right.
$$

Then $(\mathcal P)$ admits a unique solution $\eta^*=(\eta_u)_{u \in S}$.

\end{prop}

Proof of Proposition \ref{dugenou} is postponed to the Appendix. 
Notice that a similar result for the Lebesgue measure is given in \cite{hooker}. 
Its purpose was to provide diagnostics for high-dimensional functions. 
Here, we will no more exploit this idea. This will be done in a forthcoming work. 
Instead, we are going to construct stochastic sensitivity indices based on the new decomposition (\ref{equ2}) and focus on a specific estimation method for IPDV models.


\subsection{Generalized sensitivity indices}

As stated in Theorem \ref{theor1}, under (\ref{c1}) and (\ref{c2}), the output $Y$ of the model can be uniquely decomposed as a sum of hierarchically orthogonal terms. Thus, the global variance has a simplified decomposition into a sum of covariance terms. So, we can define generalized sensitivity indices.

\begin{defi}
 The sensitivity index $S_u$ of order $|u|$ measuring the contribution of $X_u$ into the model is given by :

\begin{equation}\label{eq10}
\eqbox{
S_u=\dfrac{V(\eta_u(X_u))+ \sum_{ u \cap v \neq u,v } \mathop{Cov} (\eta_u(X_u),\eta_v(X_v))}{V(Y)}
} 
\end{equation}

\bigskip

More specifically, the first order sensitivity index $S_i$ is given by :

\begin{equation}\label{eq11}
 \eqbox{
S_i=\dfrac{V(\eta_i(X_i))+ \sum_{\substack{v \neq \emptyset \\ i \not\in v }} \mathop{Cov} (\eta_i(X_i),\eta_v(X_v))}{V(Y)}
} 
\end{equation}

\end{defi}

An immediate consequence is given in Proposition \ref{pro2} (see proof in the Appendix)  :

\begin{prop}\label{pro2}
Under (\ref{c1}) and (\ref{c2}), the sensitivity 
indices $S_u$ previously defined are summed to $1$, i.e.

\begin{equation}\label{eq12}
 \sum_{ u \in S\setminus\{\emptyset\}} S_u=1
\end{equation}

\end{prop}

Thus, sensitivity indices are summed to $1$. Furthermore, the covariance terms included in these new indices allow to take into account the inputs dependence. Thus, we are now able to measure the influence of a variable on the model, especially when a part of its variability is embedded 
into the one of other dependent terms. We can distinguish the full contribution of a variable and its contribution into another correlated income.\\

Note that for independent inputs, the summands $\eta_u$ are mutually orthogonal, so $\mbox{Cov}(\eta_u,\eta_v)=0$, $u\neq v$, and we recover the well known Sobol indices. Hence, these new sensitivity indices can be seen as a generalization of Sobol indices.\\

However, the HOFD and subsequent indices are only obtained under constraints (\ref{c1}) and (\ref{c2}).
In the following, we give illustrations of distribution functions satisfying these main assumptions.

\section{Examples of distribution function}\label{sect3}

This section is devoted to examples of distribution function satisfying (\ref{c1}) and (\ref{c2}). The first hypothesis only implies that the reference measure is a product of measures, whereas the second is trickier to obtain.\\

In the first part, we  give a sufficient condition to get (\ref{c2}) for any number $p$ of input variables. The second part deals with the case $p=2$, for which we give equivalences of (\ref{c2}) in terms of copulas.

\subsection{Boundedness of the inputs density function}

The difficulty of Condition (\ref{c2}) is that the inequality has to be true for any splitting of the set  $(X_1,\cdots,X_p)$ into two disjoint blocks. We give a sufficient condition for (\ref{c2}) to hold in Proposition \ref{p5} (the proof is postponed to the Appendix):

\begin{prop}\label{p5}
Assume that there exist $M_1$,$M_2>0$ with
\cond{
$
 M_1\leq p_X\leq M_2
$
}\label{c3}
 Then, Condition (\ref{c2}) holds.
\end{prop}

Let give now an example where (\ref{c3}) is satisfied.

\exple{
Let $\nu$ be the multidimensional gaussian distribution $N_p(m,\Sigma)$ with

$$
\begin{array}{cc}
m=\begin{pmatrix}
   m_1 \\
\vdots\\
m_p
  \end{pmatrix}, & 
\Sigma = \begin{pmatrix}
          \sigma_1^2 & \cdots & 0  \\
	  & \ddots \\
	  0 & \cdots& \sigma_p^2
         \end{pmatrix}
\end{array}
$$ 

Assume that $P_X$ is a Gaussian mixture $\alpha\cdot N_p(m,\Sigma)+(1-\alpha)\cdot N_p(\mu,\Omega)$,  $\alpha \in ]0,1[$ with

$$
\begin{array}{cc}
\mu=\begin{pmatrix}
   \mu_1 \\
\vdots\\
\mu_p
  \end{pmatrix}, & 
\Omega = \begin{pmatrix}
          \varphi_1^2 & \rho_{12}&\cdots & \rho_{1p} \\
	  & \cdots \\
	  \rho_{1p} & \cdots& \cdots &\varphi_p^2
         \end{pmatrix}
\end{array}
$$ 

Then, (\ref{c3}) holds iff the matrix $(\Omega^{-1}-\Sigma^{-1})$ is positive definite.
}\label{pro5}\\

In the next section, we will see that (\ref{c2}) has a copula version when $p=2$. We will give some examples of distribution satisfying one of these conditions.

\subsection{Examples of distribution of two inputs}

Here, we consider the simpler case of inputs $\mathbf X=(X_1,X_2)$.  Also, until Section \ref{sect4}, we will assume that $\nu$ is absolutely continuous with respect to Lebesgue measure. 
The structure of dependence of $X_1$ and $X_2$ can be modelized by copulas.
Copulas~\cite{nelsen} give a relationship between a joint distribution and its marginals.
 Sklar's theorem~\cite{sklar} ensures that for any distribution function $F(x_1,x_2)$ with marginal distributions $F_1(x_1)$ and $F_2(x_2)$, $F$ has the copula representation, 
\[
 F(x_1,x_2)=C(F_1(x_1),F_2(x_2))
\]
where the measurable function $C$ is unique whenever $F_1$ and $F_2$ are absolutely continuous.\\

The next corollary gives in the absolutely continuous case the relationship between a joint density and its marginal:

\begin{cor}
In terms of copulas, the joint density of $\mathbf X$ is given by:

 \begin{equation}\label{copula1}
 p_{X}(x_1,x_2)=c(F_1(x_1),F_2(x_2))p_{X_1}(x_1)p_{X_2}(x_2)
\end{equation}

Furthermore,

\begin{equation}\label{copula2}
  c(u,v)=\dfrac{\partial^2 C}{\partial u \partial v}(u,v), \quad (u,v) \in [0,1]^2
\end{equation}
\end{cor}

Now, Condition (\ref{c2}) may be rephrased in terms of copulas:

\begin{prop}\label{copcond3}
For a two-dimensional model, the three following conditions are equivalent:
\begin{enumerate}
\item \cond{
$ p_{X} \geq M\cdot p_{X_1}p_{X_2} \quad \nu\textrm{-a.e.} \quad \textrm{for some } 0<M<1$
}\label{condmix}

\item \cond{
$
 c(u,v)\geq M, \quad \forall~(u,v)\in [0,1]^2
$
}\label{c4}
\item \cond{
$
C(u,v)=Muv+(1-M) \tilde C(u,v),
$
for some copula $\tilde C$}\label{condcop}
\end{enumerate}

\end{prop}
The proof of Proposition \ref{copcond3} is postponed to the Appendix.

Hence, the generalized Hoeffding decomposition holds for a wide class of examples. 
The first example is the Morgenstern copulas~\cite{morgenstern}:

\exple{ 
The expression of the Morgenstern copulas is given by:
\[
 C_\theta(u,v)=uv[1+\theta(1-u)(1-v)], \quad \theta \in [-1,1]
\]
For $\theta \in ]-1,1[$, (\ref{condcop}) holds, and 
\[
C_\theta(u,v)=(1-|\theta|) uv+ |\theta| uv[1+ \frac{\theta}{|\theta|} (1-u)(1-v)] 
\]
}

Let now consider the class of Archimedian copulas,

\begin{equation}\label{cop}
 C(u,v)= \varphi^{-1}[\varphi(u)+\varphi(v)], \quad u, v \in [0, 1]
\end{equation}
where the generator $\varphi$ is a non negative two times differentiable function defined on $[0,1]$ with $\varphi(1)=0$, $\varphi'(u)<0$ and  $\varphi''(u)>0,~ \forall~u \in [0,1]$.\\

A sufficient condition for (\ref{c4}) is given in Proposition \ref{procopula1}:

\begin{prop}\label{procopula1}
 If there exist $M_1$, $M_2 >0$ such that:
\begin{eqnarray}
 -\varphi'(u)\geq M_1  ~\forall~u \in [0,1]\label{con1}\\
\dfrac{d}{du}(\dfrac{1}{2}\dfrac{1}{\varphi'(u)^2})\geq M_2 , ~\forall~u \in [0,1]\label{con2}
\end{eqnarray}
 Then, Condition (\ref{c4}) holds.
\end{prop}

The proof is straightforward. Now, we will see three illustrative Archimedian copulas satisfying (\ref{c4}).

\exple{
The Frank copula is characterized by the generator:
\[
 \varphi_1(x)=\log\left( \frac{e^{-\theta x}-1}{e^{-\theta}-1}\right), \quad \theta \in \mathbb R\setminus \{0\}
\]
 Condition (\ref{c4}) holds, and $c(u,v) \geq -\theta(e^{-\theta}-1)e^{-2\theta}$ if $\theta>0$, $c(u,v)\geq -\theta(e^{-\theta}-1)$ elsewhere.
}\\

The next two examples also satisfy (\ref{c4}) by the intermediate Proposition \ref{procopula1}.

\exple{
 Let $ \alpha<0$, $\theta>0$ and $\beta$ with $\beta<-\alpha e^{-\theta}$. Set
\begin{equation}
 \varphi_2(x)=-\dfrac{\alpha}{\theta}e^{-\theta x}+\beta x+(\dfrac{\alpha}{\theta}e^{-\theta}-\beta), \quad x \in [0,1]
\end{equation}
}

\exple{
 Let  $C<0$ and set
\begin{equation}
 \varphi_3(x)=x\ln x+(C-1)x+(1-C), \quad x \in [0,1]
\end{equation}
}

Leaving the class of copulas, we now directly work with the joint density function. Proposition \ref{mixpro} gives a general form of distribution for our framework:

\begin{prop}\label{mixpro}
If $p_X$ has the form
\begin{equation}
p_{X}(x_1,x_2)=\alpha \cdot f_{X_1}(x_1)f_{X_2}(x_2)+(1-\alpha)\cdot g_{X}(x_1,x_2), \quad \alpha \in ]0,1[
\end{equation}
 where $f_{X_1}$, $f_{X_2}$ are univariate density functions, and $g_X$ is any density function (with respect to $\nu$) with marginals $f_{X_1}$ and $f_{X_2}$, then $p_X$ satisfies (\ref{c4}).
\end{prop}

 The proof is straightforward. 
\exple{ As an illustration of Proposition \ref{mixpro}, take $\nu=\nu_L$, $f_X=f_{X_1}f_{X_2}$ a normal density with a diagonal covariance matrix $\Sigma$, and $g_X$ a normal density of covariance matrix $\Omega$, with $\Omega_{ii}=\Sigma_{ii}$, $i=1,2$. 
Notice that because a copula of Gaussian mixture distribution is a mixture of Gaussian copulas (see~\cite{ouyang}), this example can be directly recovered by the copula approach.\\
}\label{exple6}

\exple{ Let generalize Example \ref{exple6}. If $P_X$ is a Gaussian mixture
\[
\alpha\cdot N_2(m,\Sigma)+(1-\alpha)\cdot N_2(\mu,\Omega), \quad \alpha \in ]0,1[
\] 

with 
 $$
\begin{array}{llll}
m=\begin{pmatrix}
   m_1 \\
m_2
  \end{pmatrix}, &
\mu=\begin{pmatrix}
   \mu_1 \\
\mu_2
  \end{pmatrix}, &
\Sigma = \begin{pmatrix}
          \sigma_1^2 &  0  \\
	  0 &  \sigma_2^2
         \end{pmatrix},  &
\Omega = \begin{pmatrix}
          \omega_1^2 & \rho \omega_1\omega_2 \\
	\rho \omega_1\omega_2  & \omega_2^2 
         \end{pmatrix},~ \rho \in ]-1,1[
\end{array}
$$
 
then (\ref{condmix}) holds iff $\omega_1^ 2 \leq \sigma_1^2$ and $\omega_2^ 2 \leq \sigma_2^2$. \\
 }

Thus, for many distributions, the generalized decomposition holds, and generalized sensitivity indices may thus be defined.\\ 

For the remaining part of the paper, we will assume that the set of inputs is an IPDV. 
If $p$ is odd, we will assume that an input variable is independent to all the others.\\
The next section is devoted to the estimation of HOFD components. The simplest case of a single pair of dependent variables is first discussed.
Then, the more general IPDV case is studied. 
In this last part, first and second order indices are defined to measure the contribution of each pair of dependent variables and each of its components in the model. 
Indices of order greater than one involving variables from different pairs will not not be studied here.

\section{Estimation}\label{sect4}
Using the property of hierarchical orthogonality ($H_u^0 \perp H_v^0$, $\forall~v \subset u $), we will see that the summands of the decomposition are solution of a functional linear system. 
For $u \in S$, the projection operator onto $H_u^0$ is denoted by $P_{H_u^0}$.\\

In this section, we present the HOFD terms computation, based on the resolution of a functional linear system. 
The result relies on projection operators previously set up. 
Further, we expose the linear system estimation for practice.

\subsection{Models of $p=2$ input variables}\label{subsect3}

This part is devoted to the simple case of bidimensional models. Let
\[
 Y=\eta(X_1,X_2)
\]

Assuming that Conditions (\ref{c1}) and (\ref{c2}) both hold, we proceed as follows:

\begin{procedure}

\begin{enumerate}\label{proc1}
\item \label{oneproc} HOFD of the output:

\begin{equation}
Y=\eta_0+\eta_1(X_1)+\eta_2(X_2)+\eta_{12}(X_1,X_2)
\end{equation}

\item \label{twoproc} Projection of $Y=\eta(\mathbf X)$ on $H_u^0$, $\forall~u \subseteq \{1,2\}$. As  $H_u^0 \perp H_v^0$, $\forall~v \subset u $, we obtain:

\begin{equation}\label{syst}
\begin{pmatrix}
\id{Id} & 0 & 0 & 0\\
 0 & \id{Id} & P_{H_1^0} & 0\\
0 & P_{H_2^0} & \id{Id} & 0\\
0 & 0 &0 & \id{Id}
\end{pmatrix}
\begin{pmatrix}
\eta_0\\
 \eta_1\\
\eta_2\\
\eta_{12}
\end{pmatrix}
=
\begin{pmatrix}
P_{H_{\emptyset}}(\eta)\\
P_{H_1^0}(\eta)\\
 P_{H_2^0}(\eta)\\
P_{H_{12}^0}(\eta)
\end{pmatrix}
\end{equation}

\item Computation of the right-hand side vector of (\ref{syst}):
\begin{equation}\label{eqs}
\begin{pmatrix}
P_{H_{\emptyset}}(\eta)\\
P_{H_1^0}(\eta)\\
 P_{H_2^0}(\eta)\\
P_{H_{12}^0}(\eta)
\end{pmatrix}=
 \begin{pmatrix}
\mathbb E(\eta)\\
\mathbb{E}(\eta/X_1)-\mathbb{E}(\eta)\\
\mathbb{E}(\eta/X_2)-\mathbb{E}(\eta)\\
\eta-\mathbb{E}(\eta/X_1)-\mathbb{E}(\eta/X_2)+\mathbb{E}(\eta)
\end{pmatrix}
\end{equation}

In this frame, we have:

\begin{prop}\label{pro3}
 Let $\eta$ be any function of $L^2_{\mathbb R}(\mathbb R^p,\mathcal{B}(\mathbb{R}^p),P_X)$. 
Then, under (\ref{c1}) and (\ref{c2}), the linear system  

\begin{equation}\label{eq15}
(\mathcal S)
\begin{pmatrix}
\id{Id} & 0 & 0 & 0\\
 0 & \id{Id} & P_{H_1^0} & 0\\
0 & P_{H_2^0} & \id{Id} & 0\\
0 & 0 &0 & \id{Id}
\end{pmatrix}
\begin{pmatrix}
h_0\\
 h_1\\
h_2\\
h_{12}
\end{pmatrix}
=
\begin{pmatrix}
P_{H_{\emptyset}}(\eta)\\
P_{H_1^0}(\eta)\\
 P_{H_2^0}(\eta)\\
P_{H_{12}^0}(\eta)
\end{pmatrix}
\end{equation}

admits in $h=(h_0,\cdots,h_{12}) \in H_{\emptyset}\times \cdots \times H_{12}^0$ the unique solution $h^*=(\eta_0,\eta_1,\eta_2,\eta_{12})$.

\end{prop}

\item Reduction of the system (\ref{syst}). As the constant term corresponds to the expected value of $\eta$, and the residual one can be deduced from the others, 
the dimension of the system (\ref{eq15}) can even be reduced to:

\begin{equation}\label{eq16}
\underbrace{
\begin{pmatrix}
\id{Id} & P_{H_1^0} \\
P_{H_2^0} & \id{Id}\\
\end{pmatrix}
}_{A_2}
\underbrace{
\begin{pmatrix}
 \eta_1\\
\eta_2\\
\end{pmatrix}
}_{\Delta}
=
\underbrace{
\begin{pmatrix}
\mathbb{E}(\eta/X_1)-\mathbb{E}(\eta)\\
\mathbb{E}(\eta/X_2)-\mathbb{E}(\eta)
\end{pmatrix}
}_{B}
\end{equation}

\item Practical resolution: The numerical resolution of (\ref{eq16}) is achieved by an iterative Gauss Seidel algorithm~\cite{gs} which consists first in
decomposing $A_2$ as a sum of lower triangular ($L_2$) and strictly upper triangular ($U_2$) matrices.

Further, the technique uses an iterative scheme to compute $\Delta$. At step $k+1$, we have:

\begin{equation}
 \Delta^{(k+1)}: =
\begin{pmatrix}
 \Delta_1^{(k+1)} \\
 \Delta_2^{(k+1)}
\end{pmatrix} = L_2^{-1}(B-U_2\cdot\Delta^{(k)})
\end{equation}
Using expression of $A_2$, we get: 

\begin{equation}\label{scheme}
   \Delta^{(k+1)}=
 \begin{pmatrix}
\mathbb E(Y-\Delta_2^{(k)}/X_1)-\mathbb E(Y-\Delta_2^{(k)}) \\
\mathbb E(Y- \Delta_1^{(k+1)}/X_2)-\mathbb E(Y- \Delta_1^{(k+1)})
\end{pmatrix}
\end{equation}

\item \label{sixproc} Convergence of the algorithm: now, we hope that the Gauss Seidel algorithm converges to the true solution. Looking back at (\ref{syst}) , we see that we only have to consider $P_{H_1^0}$ (respectively $P_{H_2^0}$) restricted to $H_2^0$ (respectively to $H_1^0$).

Under this restriction, let us define the associated norm operator as :

\[
 \| P_{H_i^0}\|^2 := \sup_{\substack{\mathbb E(U^2)=1 \\ U \in H_j^0}} \mathbb E[P_{H_i^0}(U)^2], \quad i,j=1,2, ~j \neq i
\]

As explained in \cite{pillis}, Gauss Seidel algorithm converges to the true solution $\Delta$ if  $A_2$ is striclty diagonally dominant, which is implied by :

\begin{equation}\label{n1}
 \| P_{H_i^0}\|<\|\id{Id}\|, \quad i=1,2
\end{equation}

i.e.

\begin{equation}\label{n2}
 \sup_{\substack{\mathbb E(U^2)=1 \\ U \in H_j^0}} \mathbb E[P_{H_i^0}(U)^2] < 1, \quad i=1,2, ~j \neq i
\end{equation}

Equality (\ref{eqs}) reveals that $ P_{H_i^0}(U)=\mathbb E(U/X_i)-\mathbb E(U)$. 
Hence, by the Jensen inequality~\cite{jen} for conditional expectations, $ \| P_{H_i^0}\|, i=1,2$ admits an upper bound:\\

Take $U \in H_1^0$ :

$$
\begin{array}{lll}
 \| P_{H_2^0}\|&= &\displaystyle{\sup_{\substack{\mathbb E(U^2)=1 \\ U \in H_1^0}} \mathbb E[(\mathbb E(U/X_2)-\mathbb E(U))^2]} \\
&=&\displaystyle{\sup_{\substack{\mathbb E(U^2)=1 \\ U \in H_1^0}} \mathbb E[\mathbb E(U/X_2)^2]\quad \textrm{ as } U \in H_1^0}\\
&\leq & \displaystyle{\sup_{\substack{\mathbb E(U^2)=1 \\ U \in H_1^0}} \mathbb E[\mathbb E(U^2/X_2)]}=1 \quad \textrm{by Jensen}\\
\end{array}
$$
 The same holds for $U \in H_2^0$, and we also have $ \| P_{H_2^0}\| \leq 1$.\\
 Moreover, the Jensen's inequality is strict if $U$ is not $X_i$-measurable. 
As $U$ is a function of $X_j$ (that is $j=2$ if $i=1$ and conversely), 
the condition of convergence holds if $X_1$ is not a measurable function of $X_2$.
Hence, Gauss Seidel algorithm converges if $X_1$ is not a function of $X_2$.\\

\item \label{numproc} Estimation procedure: 
Suppose that we get a sample of $n$ observations $(Y_k,\mathbf X_k)_{k=1,\cdots,n}$. 
\begin{itemize}
\item Estimation of the components of the HOFD: 
the iterative scheme (\ref{scheme}) requires to estimate conditional expectations.
As extended in Da Veiga {\it{et al}}~\cite{daveiga}, we propose to estimate them by local polynomial regression at each point of observation.
 Then, we use the leave-one-out technique to set the learning sample and the test sample. Moreover, as the local polynomial method can be summed up to a generalized least squares 
(see Fan and Gijbels~\cite{fan}), the Sherman-Morrison formula~\cite{SM} is applied to reduce the computational time.\\
A more detailed procedure is given in the Appendix. 
The iterative algorithm is easy to implement. We stop when $\| \Delta^{(k+1)}- \Delta^{(k)} \|\leq \varepsilon$, for a small positive $\varepsilon$.

Once $(\eta_1,\eta_2)$ have been estimated, we estimate $\eta_0$ by the empirical mean of the output. Then, an estimation of $\eta_{12}$ is obtained by substraction.

\item We use empirical variance and covariance estimation to estimate sensitivity indices $S_1$, $S_2$ and $S_{12}$.
\end{itemize}

\end{enumerate}

\end{procedure}


\subsection{Generalized IPDV models}

Assume that the number of inputs is even, so $p=2k$, $k\geq 2$. 
We note each group of dependent variables as $X^{(i)}:=(X_1^{(i)},X_2^{(i)})$, $i=1,\cdots,k$.
By rearrangement, we may assume that:

\[
 \mathbf X=(\underbrace{X_1,X_2}_{X^{(1)}}, \cdots, \underbrace{X_{2k-1},X_{2k}}_{X^{(k)}})
\]

SA for IPDV models has already been treated in \cite{jac}. 
Indeed, they proposed to estimate usual sensitivity indices on groups of variables via Monte Carlo estimation.
Thus, they have interpreted the influence of every group of variables on the global variance.
Here, we will go further by trying to measure the influence of each variable on the output, but also the effets of the independent pairs.\\
To begin with, as a slight generalization of \cite{sobol} and used in \cite{jac}, let apply the Sobol decomposition on groups of dependent variables,
\[
  \eta(\mathbf{X})=\eta_0+\eta_1(X^{(1)})+\cdots+\eta_k(X^{(k)})+\sum_{|u|= 2}^k \eta_u(X^{(u)})
\]

where for $u=\{u_1,\cdots,u_t\}$ and $t=|u|$, we set $X^{(u)}=(X^{(u_1)},\cdots,X^{(u_t)})$. 
Furthermore, $\langle \eta_u,\eta_v \rangle=0$, $\forall~u \neq v$.\\

Thus, we obtain independent components of IPDV. Under the assumptions discussed in the previous section, we can apply the HOFD of each of these components, that is,
\[
\eta_i(X^{(i)})=\eta_i(X^{(i)}_{1},X^{(i)}_{2})=\varphi_{i0}+\varphi_{i,1}(X^{(i)}_{1})+\varphi_{i,2}(X^{(i)}_{2})+\varphi_{i,12}(X^{(i)})
\]

with $\langle \varphi_{i,u},\varphi_{i,v}\rangle=0$, $\forall~v \subset u \subseteq \{1,2\}$.
In this way, let define some new generalized indices for IPDV models:

\begin{defi}
 For $i=1,\cdots,k$, the first order sensitivity index measuring the contribution of $X_1^{(i)}$ (respectively $X_2^{(i)}$) on the output of the model is:

\begin{equation}\label{ipdvs1}
  S_{i,1}=\dfrac{V(\varphi_{i,1})+\mbox{Cov}(\varphi_{i,1},\varphi_{i,2})}{V(Y)}, ~ \left(S_{i,2}=\dfrac{V(\varphi_{i,2})+\mbox{Cov}(\varphi_{i,1},\varphi_{i,2})}{V(Y)}\right)
\end{equation}

The second order sensitivity index for the pair $X^{(i)}$, $i=1,\cdots,k$, is defined as:

\begin{equation}\label{ipdvs2}
 S_{i,12}=\dfrac{V(\varphi_{i,12})}{V(Y)}
\end{equation}

\end{defi}

The estimation procedure of these indices is quite similar to Procedure \ref{proc1}:

\begin{procedure}

\begin{enumerate}\label{proc2}

\item Estimation of $(\eta_i)_{i=1,\cdots,k}$: as reminded in Part \ref{sect22} with Equations (\ref{sd}), $\eta_i=\mathbb E(Y/X^{(i)})-\mathbb E(Y)$. Step \ref{numproc} of Procedure \ref{proc1} gives method to estimate the conditional expectations.  So that, we will have estimations of $\eta_i$, $i=1,\cdots,k$.

\item  For $i=1,\cdots,k$, we apply step \ref{twoproc} to step \ref{numproc} of Procedure \ref{proc1}, considering $\eta_i$ as the output.

\end{enumerate}

\end{procedure}
 If $p$ is odd, the procedure is the same except that the influence of the independent variable is measured by a Sobol index, as it is independent from all the others.\\
The next part is devoted to numerical examples.

\section{Numerical examples}\label{sect5}
In this section, we study three examples with dependent input variables. 
We consider IPDV models and a Gaussian mixture distribution on the input variables. 
We choose covariance matrices of the mixture satisfying conditions of Example \ref{pro5}.\\
We give estimations of our new indices, and compare them to the analytical ones, computed from expressions (\ref{eq10}). We also compute dispersions of the estimated new indices.\\
In~\cite{daveiga}, Da Veiga {\it{et al.}} proposed to estimate the classical Sobol indices $ S_u=(V[\mathbb E(Y/X_u)]-\sum_{v \subset u} V[\mathbb E(Y/X_v)])/V(Y), u \subseteq \mathcal P_p$, by nonparametric tools.
Indeed, the local polynomial regression were used to estimate conditional moments $\mathbb E(Y/X_u), u \subseteq \mathcal P_p$. This method, used further, will be called Da Veiga procedure (DVP). 
Results given by DVP are compared with the ones given by our method. 
The goal is to show that the usual sensitivity indices are not appropriate in the dependence frame, even if a relevant estimation method is used.

\subsection{Two-dimensional IPDV model}

Let consider the model

\[
Y=X_1+X_2+X_1X_2
\]

Here, $\nu$ and $P_X$ are of the form given by Example \ref{pro5}, with $m=\mu=0$. \\

Thus, the analytical decomposition of $Y$ is

\begin{equation*}
  \eta_0=\mathbb E(X_1X_2),\quad \eta_1=X_1,\quad \eta_2=X_2 \quad \eta_{12}=X_1X_2 -\mathbb E(X_1X_2)
\end{equation*}

For the application, we implement Procedure \ref{proc1} in Matlab software. 
We proceed to $L=50$ simulations and $n=1000$ observations. 
Parameters were fixed at $\sigma_1=\sigma_2=1$, $\varphi_1^2=\varphi_2^2=0.5$, $\rho_{12}=0.4$ and $\alpha=0.2$.\\

In Table \ref{tab1}, we give the estimation of our indices and their standard deviation (indicated by $\pm\cdot$) on $L$ simulations.
In comparison, we give the analytical value of each index. 

The analytical classical Sobol indices are difficult to obtain, but we give estimators of the classical Sobol indices with DVP.

\begin{table}[H]
\caption{Estimation of the new and DVP indices with $\rho_{12}=0.4$}\label{tab1}
\centering
\renewcommand{\arraystretch}{1.9}
\begin{tabular}{|l|c|c|c|c||c|}
 \hline
& & $S_1$ & $S_2$ & $S_{12}$ & $\sum_u S_u$\\
\hline
\multirow{2}{1cm}{New indices} & Estimation & $0.42 \pm 0.041$ &  $0.41 \pm 0.043$   & $0.17\pm 0.026$ & $1\pm 9.10^{-16}$\\
\cline{2-6}
& Analytical & $0.39 $ & $0.39$ & $0.22$ & $1$\\
\hline\hline 
\multirow{1}{1cm}{DVP indices} &   Estimation & $0.64 \pm 0.045$  &   $0.65\pm 0.044$  &    $0.41 \pm 0.038$  & $1.7\pm 0.09$\\
\hline
\end{tabular}
\end{table}

We notice that estimations with our method give quite good results in comparison with their analytical values. 
The estimation error of the interaction term is due to the fact that the component $\hat\eta_{12}$ is obtained by difference between the output and the other estimated components.

The DVP indices are are difficult to interpret as the sum is higher than $1$.

In our method, it would be relevant to separate the variance part to the covariance one in the first order indices. 
Indeed, in this way, we would be able to get the part of variability explained by $X_i$ alone in $S_i$, and its contribution hidden in the dependence with $X_j$.
We note $S_i^v$ the variance contribution alone, and $S_i^{c}$ the covariance contribution, that is

\[
 S_i=\underbrace{\dfrac{V(X_i)}{V(Y)}}_{S_i^v}+\underbrace{\dfrac{\mbox{Cov}(X_i,X_j)}{V(Y)}}_{S_i^c}, \quad i=1,2, ~j\neq i
\]

The new indices estimations given in Table \ref{tab1} are decomposed in Table \ref{tab11}. 
As previously, the number at the right of $\pm$ indicates the standard deviation on $L$ simulations.

\begin{table}[H]
\caption{Estimation of $S_i^v$ and $S_i^c$ with $\rho_{12}=0.4$}\label{tab11}
\centering
\renewcommand{\arraystretch}{1.5}
\begin{tabular}{|l|c|c|c|}
 \hline
& $S_i^v$ & $S_i^c$ & $S_i$ \\
\hline
$X_1$ &  $ 0.28 \pm 0.04$ &  $0.14\pm 0.01 $   & $0.42 \pm 0.041$ \\
\hline
$X_2$ &  $0.27 \pm 0.043$ & $0.14\pm 0.01$ & $0.41 \pm 0.043$\\
\hline\hline 
Analytical & $0.25 $ & $0.14$ & $0.39$\\
\hline
\end{tabular}
\end{table}

For each index, the covariate itself explains $28\%$ (in estimation, $25\%$ in reality) of the part of the total variability. 
However, the contribution embedded in the correlation is not negligible as it represents $14\%$ of the total variance.
Considering the shape of the model, and coefficients of parameters distribution, it is quite natural to get the same contribution of $X_1$ and $X_2$ into the global variance. 
Also, as their dependence is quite important with a covariance term equals to $0.4$, we are not surprised by the relatively high value of $S_1^c$ (resp. $S_2^c$). \\

From now, we take $\rho_{12}=0$, i.e. we assume that the inputs are independent. 
Let compare our new estimated indices with their analytical values in Table \ref{tab2}. 
We again decompose new indices into a variance ($S_i^v$) and a covariance ($S_i^c$) contribution.

\begin{table}[H]
\caption{Comparison between analytical and estimated indices with $\rho_{12}=0$}\label{tab2}
\centering
\renewcommand{\arraystretch}{1.5}
\begin{tabular}{c|c|c|c||c|}
 \cline{2-5}
& \multicolumn{3}{c||}{New indices}  & \multicolumn{1}{c|}{Theoretical}\\
\cline{2-5}
& $S_i^v$ & $S_i^c$ & $S_i$ & $S_i$ \\
\hline
\multicolumn{1}{|c|}{$X_1$} &  $ 0.39 \pm 0.039$ &  $-0.01\pm 0.01 $   & $0.38 \pm  0.036$ & $0.375$\\
\hline
\multicolumn{1}{|c|}{$X_2$} &  $0.38 \pm 0.045$ & $-0.01\pm 0.024$ & $0.37 \pm  0.037 $& $0.375$\\
\hline
\multicolumn{1}{|c|}{$X_{1}X_{2}$} & $0.26\pm 0.038$ & $-0.01\pm 0.01$ &  $0.25\pm 0.027$ & $0.25$ \\
\hline
\end{tabular}
\end{table}

Thus, the new indices are well tailored if we have a small idea on inputs dependence in a system. 
Indeed, Table \ref{tab2} shows that our new indices take dependence into account if it exists, and the covariance contribution is estimated by $0$ if not. 
New indices recover the classical Sobol indices in case of independence.


\subsection{Linear Four-dimensional model}

The test model is

\[
Y=5X_1+4X_2+3X_3+2X_4
\]

Actually, Condition (\ref{c2}) only needs to be satisfied on groups of correlated variables. \\
Let consider the two blocks $X^{(1)}=(X_1,X_3)$ and $X^{(2)}=(X_2,X_4)$ of correlated variables. \\
The previous form of density can be taken for $X^{(1)}$ and $X^{(2)}$. 
$P_{X^{(i)}}$ is then the Gaussian mixture $\alpha_i\cdot N_2(0,\id{I}_2)+(1-\alpha_i)\cdot N_2(0,\Omega_i)$, $i=1,2$.
The analytical sensitivity indices are given by (\ref{ipdvs1}) \& (\ref{ipdvs2}).\\

For $L=50$ simulations and $n=1000$ observations, we took $\varphi_{1}^{2(1)}=\varphi_{2}^{2(1)}=0.5$, $\varphi^{2(2)}_{1}=0.7$, $\varphi^{2(2)}_{2}=0.3$, $\rho_{12}^{(1)}=0.4$, $\rho_{12}^{(2)}=0.37$ and $\alpha_1=\alpha_2=0.2$.\\

\begin{figure}[H]
\centering
 \includegraphics[width=0.75\textwidth,angle=-90]{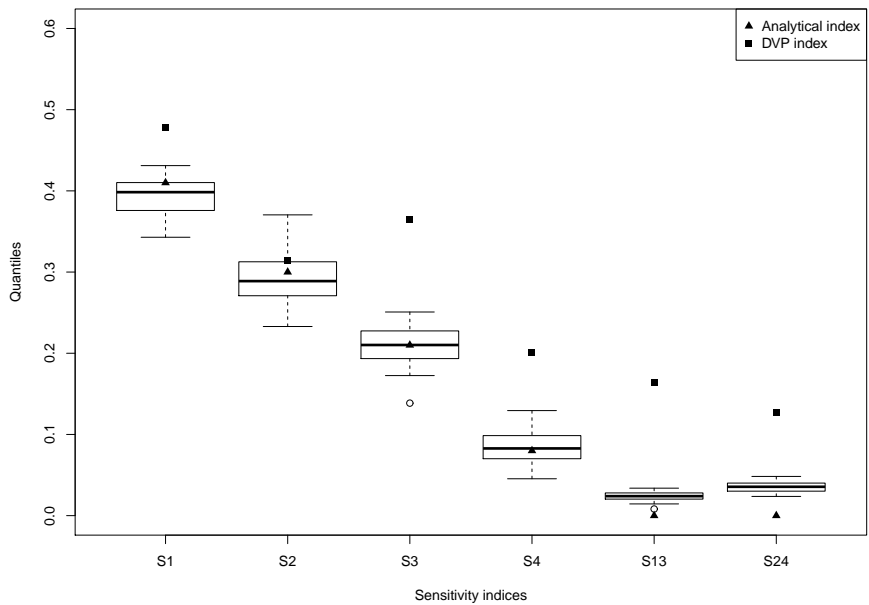}
\caption{Boxplots representation of new indices-Comparison with analytical and DVP indices}\label{fig2}
\end{figure}

Figure \ref{fig2} displays the dispersion of indices of first order for all variables and second order for grouped variables. We compare them to their analytical values.
In the same figure, we also represented the estimators of classical Sobol indices with DVP.

 We see that $X_1$ has the biggest contribution, whereas the influence of $X_4$ is very low. 
It reflects well the model if we look at the coefficients of $X_i$, $i=1,\cdots,4$. 
Also, interaction terms are well estimated, as they are closed to $0$.
For each case, the dispersion on $50$ simulations is very low.\\
As for the DVP estimation, it is once again very high compared with the true indices values.

\subsection{The Ishigami function}

This function is well known in SA (\cite{ishigami}). It is defined by:

\[
 Y=\sin(X_1)+a\sin^2(X_2)+bX_3^3\sin(X_1)
\]
We assume that $X_3$ is the independent variable, and that $X_1$ and $X_2$ are correlated. 
$P_X$ is again the Gaussian mixture $\alpha\cdot N_3(0,I_3)+(1-\alpha)\cdot N_3(0,\Omega)$.

\begin{figure}[H]
\centering
 \includegraphics[width=0.8\textwidth,angle=-90]{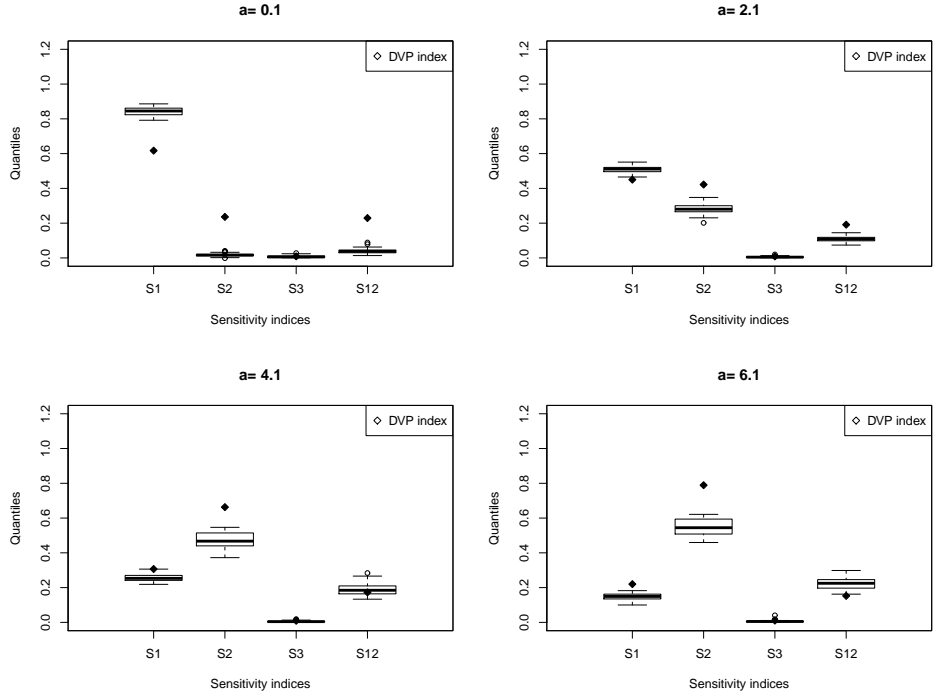}
\caption{Boxplots of new indices for different values of $a$ for Ishigami function}\label{fig1}
\end{figure}

With $L=50$ simulations of $n=1000$ observations, we fixed parameters of distribution at $\varphi_1^2=0.15$, $\varphi_2^2=0.85$, $\varphi_3^2=0.75$, $\rho_{12}=0.3$ and $\alpha=0.2$.\\
In Figure \ref{fig1}, the dispersion of the new measures is represented for fixed $b=0.1$ and different values of $a$.
In addition, the analytical new indices and estimated classical Sobol indices with DVP are displayed.

Every boxplot shows that there is a small dispersion. 
For all estimations, DVP indices are larger that the new ones. 
The figure clearly shows that, for all values of $a$, the sum of these four indices is greater 
than $1$. It again shows their non adaptation to a situation of dependence.\\
If we have a look on the values taken by our new sensitivity indices, we see that, 
for small values of $a$, the variable $X_1$ contributes the most to the model's variability. 
This role decreases as $a$ increases, and $X_2$ then gets the biggest contribution. 
For any value of $a$, the input $X_3$ plays a very negligible role, which seems realistic as $b$ is a small fixed value.
As for the interaction index $S_{12}$, it is getting bigger with the increasing importance of $a$, but the contribution remains low.


\section{Conclusions and Perspectives}\label{sect6}

The Hoeffding decomposition and associated Sobol indices have been widely studied in SA over past years. 
Recently, a literature appears to treat the case of dependent input variables, in which authors propose different ways to deal with dependence. 
The goal of this paper is to conciliate the problem of inputs dependence with the Hoeffding decomposition. 
Indeed, we study a functional ANOVA decomposition in a generalized inputs distribution frame. 
Thus, we show that a model can be uniquely decomposed as a sum of hierarchically orthogonal functions of increasing dimension. 
Also, this approach generalizes the Hoeffding's one, as we recover it in case of independence.\\ 
Similarly to the classical Sobol decomposition, this leads to the construction of new sensitivity indices. 
They consist of a variance and a covariance contribution able to take into account the possible correlation among variables.
In case of independence, these indices are the classical Sobol indices. 
However, the indices construction is only possible under specific assumptions on the joint distribution function of the inputs. 
We expose few cases that satisfy these assumptions for any $p$-dimensional models. 
More specifically, for two-dimensional models, the required assumption is equivalent to assumptions on copulas.
In this context, we give examples satisfying one of these assumptions .\\
Focused on the IPDV models, summands of the decomposition are estimated thanks to projection operations. 
This leads to the numerical resolution of functional linear systems. 
The strength of this method is that it does not require to make assumptions on the form of the model or on the structure of dependence. We neither use meta-modelling and avoid in this way many sources of errors.\\
Through three applications on test models, we observe the importance of considering the inputs correlation, and show how our method could catch it. 
The comparison with estimators of classical indices with DVP shows that the Sobol indices are not appropriate in case of correlations, even when using nonparametric method.
Also, when inputs independence holds, the new indices remain well suited to measure sensitivity into a model.\\
Nevertheless, only considering IPDV models for estimation is restrictive. 
The perspective is to explore other estimation methods suitable for more general models. 
Also, we intend to lead a systematic study on copulas satisfying or not our assumptions.


\appendix{

\section{Generalized Hoeffding decomposition}
\subsection{Generalized decomposition for dependent inputs}\label{ap1}

\itshape{The upcoming proof follows the guideline of the proof of Lemma  3.1 in Stone~\cite{stone}.}\\

\begin{preuve}{Proof of Lemma \ref{lemstone}}

By induction on the cardinal of $T$, let show that 

$$
\begin{array}{ccc}
  \mathcal H(n) :& \forall~T/ \#(T)=n,& \mathbb E[(\sum_{u \in T }h_u(\mathbf X))^2] \geq \delta^{\#(T)-1} \sum_{u \in T} \mathbb E[h_u^2(\mathbf X)] 
\end{array}
$$

\begin{itemize}
 \item $\mathcal H(1)$ is obviously true, as $T$ is reduced to a singleton 

\item Let $n \in \mathbb N^*$. Suppose that $\mathcal H(n')$ is true for all $1\leq n' \leq n$. Let $T$ such that $\#(T)=n+1$. 
We want to prove $\mathcal H(n+1)$.

Choose a maximal set $r$ of $T$, i.e. $r$ is not a proper subset of any set $u$ in $T$. We show first that 

\begin{eqnarray}\label{lem4}
 \mathbb E[(\sum_{u \in T} h_u(\mathbf X))^2] \geq M \cdot\mathbb{E}(h_r ^2(\mathbf X))
\end{eqnarray}

\begin{itemize}
  \item If $\#(r)=p$, by definition of $H_r^0$, we get ~$ \mathbb E[(\sum_{u \in T} h_u( \mathbf X))^2]  \geq  \mathbb{E}(h_r ^2(\mathbf X)) \geq  M  \mathbb{E}(h_r ^2(\mathbf{X}))$ as  $M \leq 1$.

\item If $1 \leq \#(r) \leq p-1$, set $\mathbf{X}=(X_1,X_2)$, where $X_1=(X_l)_{l \notin r}$ and $X_2=(X_l)_{l \in r}$. By Condition (\ref{c2}), it follows that

$$
\begin{array}{l}
 p_X \geq M \cdot p_{X_1}p_{X_2}
\end{array}
$$

As a consequence,

$$
\begin{array}{lll}
 \mathbb E[(\sum_{u \in T} h_u(\mathbf X))^2] & = & \int_{\mathcal X_1} \!\! \int_{\mathcal X_2}[h_r(x_2)+ \sum_{u \neq r} h_u(x_1,x_2)]^2 p_X \nu(dx_1,dx_2)\\
& \geq & M \int_{\mathcal X_1} \!\! \int_{\mathcal X_2}[h_r(x_2)+ \sum_{u \neq r} h_u(x_1,x_2)]^2 p_{X_1}p_{X_2}\nu_1(dx_1)\nu_2(dx_2) \\
&\geq & M \int_{\mathcal X_1} \mathbb E[(h_r(X_2)+\sum_{u \neq r} h_u(x_1,X_2))^2] p_{X_1} \nu_1(dx_1)
\end{array}
$$

By maximality of $r$ and by definition of $H_r^0$, 

\begin{itemize}
 \item If $u \subset r$, $h_u$ only depends on $X_2$ and by orthogonality, 
\[
\mathbb E(h_u(X_2)h_r(X_2))=0 
\]
\item If $u  \not\subset r$, $h_u$ depends on $X_1$ fixed at $x_1$, and $X_2^u=(X_l)_{l \in r \cap u}$, so $h_u \in H_{r\cap u}^0$, with $r\cap u \subset r$, it comes then

\[
   \mathbb E(h_u(x_1,X_2)h_r(X_2))= 0 
\]

\end{itemize}

Thus,

$$
\begin{array}{lll}
 \displaystyle{ \mathbb E[(\sum_{u \in T} h_u( \mathbf X))^2]} & \geq &  \displaystyle{M \int_{\mathcal X_1}\mathbb E(h_r^2(X_2)) p_{X_1}\nu_1(dx_1)}\\
& = & \displaystyle{M \cdot \mathbb E(h_r^2(\mathbf X)) }
\end{array}
$$

\end{itemize}

So (\ref{lem4}) holds for any size of any maximal sets of ~$T$.\\

By using (\ref{lem4}) with $\tilde{h}_r=h_r$ and $\tilde{h}_u=-\beta h_u$, $\forall~ u\neq r$, we get

\begin{eqnarray}\label{lem5}
   \mathbb E[(h_r(\mathbf X)-\beta \sum_{u \neq r} h_u(\mathbf X))^2] \geq M \mathbb{E}(h_r ^2(\mathbf X))
\end{eqnarray}

Taking $\beta=\dfrac{\mathbb E[h_r(\mathbf X) \sum_{u \neq r}h_u(\mathbf X)]}{\mathbb E[(\sum_{u \neq r}h_u(\mathbf X))^2]}$, it follows that:
$$
\begin{array}{llll}
  &\mathbb E[(h_r(\mathbf X)-\beta \sum_{u \neq r} h_u(\mathbf X))^2] & \geq &  M \mathbb{E}(h_r ^2(\mathbf X)) \\
\textrm{that is} & \mathbb E[h_r^2(\mathbf X)]-\dfrac{[\mathbb E(h_r(\mathbf X)\sum_{u \neq r}h_u(\mathbf X))]^2}{\mathbb E[(\sum_{u \neq r}h_u(\mathbf X))^2]}  & \geq &  M \mathbb{E}(h_r ^2(\mathbf X))
\end{array}
$$

Hence,

\begin{eqnarray}\label{ine}
 \mathbb [E(h_r(\mathbf X)\sum_{u \neq r}h_u(\mathbf X))]^2 \leq (1-M)\cdot \mathbb{E}(h_r ^2(\mathbf X))\cdot \mathbb E[(\sum_{u \neq r}h_u(\mathbf X))^2]
\end{eqnarray}

This implies

\begin{eqnarray}\label{ope}
 \mathbb E[(\sum_{u}h_u(\mathbf X))^2] \geq (1-\sqrt{1-M}) \left[\mathbb{E}(h_r ^2(\mathbf{X}))+ \mathbb E[(\sum_{u \neq r}h_u(\mathbf X))^2] \right]
\end{eqnarray}

Set $x=h_r(\mathbf X) \textrm{ and }  y=\sum_{u \neq r} h_s(\mathbf X)$.(\ref{ope}) is rephrased as

\begin{eqnarray}\label{trs}
 \|x+y\|^2 \geq (1 -\sqrt{1-M})\{\|x \|^2+ \| y\|^2\}
\end{eqnarray}

Further, (\ref{ine}) is $\langle x,y\rangle \geq - \sqrt{1-M} \|x \|\cdot \| y\|$. Thus,
$$
\begin{array}{lll}
 \|x\|^2+\|y\|^2 &\geq & 2 \langle x,y\rangle \\
& \geq & -\dfrac{2}{\sqrt{1-M}}\langle x,y\rangle \quad \textrm{by (\ref{trs})} 
\end{array}
$$
So $\|x+y\|^2  \geq  (1-\sqrt{1-M}) \{\|x \|^2+ \| y\|^2\}$.\\

As $\mathcal H(n)$ is supposed to be true and (\ref{ope}) holds, it follows that:

$$
\begin{array}{lll}
  \mathbb E[(\sum_{u}h_u(\mathbf X))^2] &\geq &\delta \left[\mathbb{E}(h_r ^2(\mathbf X))+ \delta^{n-1} \sum_{u \neq r} \mathbb E(h_u^2( \mathbf X)) \right] \\
& \geq & \delta^{n} \sum_{u} \mathbb E(h_u^2(\mathbf X)) \qquad \textrm{ as }~\delta \in ~]0,1]\\
& = &  \delta^{\#(T)-1} \sum_{u} \mathbb E(h_u^2(\mathbf X))
\end{array}
$$

Hence, $\mathcal H(n+1)$ holds.
\end{itemize}
We can deduce that $\mathcal H(n)$ is true for any collection $T$ of $\mathcal P_p$.
\end{preuve}


 \label{annexpro1}

\begin{preuve}{Proof of  Theorem \ref{theor1}}

Let define the vector space $K^0=\{ \sum_{u \in S^-}h_u(X_u),h_u \in H_u^0, \forall~u \in S^- \}$.\\

In the first step, we will prove that $K^0$ is a complete space to prove the existence and uniqueness of the projection of $\eta$ in $K^0$,
thanks to the projection theorem~\cite{luen}.\\
Secondly, we will show that $\eta$ is exactly equal to the decomposition into $H^0$, and finally, we will see that each term of the summand
is unique.
\begin{itemize}
 \item We show that $H_u^0$ is closed into $H_u$ (as $H_u$ is a Hilbert space).\\

Let $(h_{n,u})_n$ be a convergent sequence of $H_u^0$ with $h_{n,u}\rightarrow h_u$. As $(h_{n,u})_n \in H_u^0\subset H_u$ complete, $h_u \in H_u$. Let $v\subset u$, and $h_v \in H_v^0$ :

$$
\begin{array}{cccc}
 \langle h_u-h_{n,u},h_v\rangle & = \langle h_u,h_v\rangle &- \langle h_{n,u},h_v \rangle &\\
\downarrow & & \shortparallel& \\
0 & & 0 & \textrm{ as }H_u^0\perp H_v^0 
\end{array}
$$
Thus, $\langle h_u,h_{v}\rangle= 0$, so that $h_u \in H_u^0$. $H_u^0$ is then a complete space.\\

Let $(h_n)_n$ be a Cauchy sequence in $K^0$ and we show that each component is of Cauchy and that $h_n\rightarrow h \in K^0$.

As $h_n \in K^0$, $h_n=\sum_{u \in S^-} h_{n,u}$, $h_{n,u} \in H_u^0$. It follows that :
$$
\begin{array}{lll}
 \| h_{n}-h_{m}\|^2 &=&  \|\sum_u( h_{n,u}-h_{m,u})\|^2 \\
& \geq & \delta^{\#(S^-)-1} \sum_{u\in S^-} \| h_{n,u}-h_{m,u}\|^2 \quad\textrm{ by Inequality (\ref{lem})}
\end{array}
$$

 As $(h_{n})_n$ is a Cauchy sequence,  by the above inequality, $(h_{n,u})_n$ is also Cauchy. 
As $h_{n,u}\rightarrow h_u \in H_u^0$, we deduce that $ h_n \underset {n\to \infty}{\longrightarrow} \sum_{u \in S^-}h_u =h \in K^0$.

Thus, $K^0$ is complete. By the projection theorem, we can deduce there exists a unique element into $K^0$ such that :

\[
 \|\eta-\sum_{u \in S^-}\eta_u\|^2 \leq \|\eta-h\|^2 \quad \forall~h \in K^0
\]

\item Decomposition of $\eta$: following Hooker~\cite{hooker}, we introduce the residual term as

\[
 \eta_{\mathcal P_p}(X_1,\cdots,X_p)=\eta(X_1,\cdots,X_p)-\sum_{u \in S^-}\eta_u(X_u)
\]

By projection, $ \langle \eta -\sum_{v \in S^-} \eta_v , h_u\rangle=0$ $\forall ~u\in S^-$, $\forall ~h_u \in H_u^0$. 
Hence, $\eta(\mathbf X)=\sum_{u \in S}\eta_u(X_u)$, $\eta_u \in H_u^0$, $\forall ~u \in S$, and this decomposition is well defined.

\item Terms of the summand are unique: assume that ~$\eta=\sum_{u \in S} \eta_u=\sum_{u \in S} \widetilde \eta_u, \quad \widetilde \eta_u \in H_u^0$.\\
By Lemma \ref{lemstone}, it follows that

$$
\begin{array}{lll}
\left.
\begin{array}{l}
\sum_{u \in S}(\eta_u-\widetilde{\eta}_u)=0 \\
\| \sum_{u \in S}(\eta_u-\widetilde{\eta}_u) \|^2 \geq \delta^{\#(S)-1}\sum_{u \in S}\|\eta_u-\widetilde{\eta}_u\|^2\\
\end{array}
 \right \} 
& \Rightarrow &
\|\eta_u-\widetilde{\eta}_u\|^2=0 \quad \forall~u \in S
\end{array}
$$
\end{itemize}
 
\end{preuve}


\label{annexdugenou}
\begin{preuve}{Proof of Proposition \ref{dugenou}}

Let first prove the following equivalence :

$$
\begin{array}{c}
  \int \eta_u(x_u)\eta_v(x_v)p_X(x)d\nu(x) =0 \quad \forall v \subset u, ~\forall ~\eta_v\\
\Updownarrow\\
 \int \eta_u(x_u)p_X(x) d\nu_i(x_i)~ d\nu_{u^c}(x_{u^c})=0 \quad \forall u \in S, ~\forall i \in u
\end{array}
$$

\bigskip
%

Let $v \subset u$ and $i \in u\setminus v$, then 

$$
\begin{array}{lll}
\int \eta_v(x_v)\eta_u(x_u) p_X(x)d\nu(x)&=& \int   \eta_v(x_v)\eta_u(x_u) p_X(x) d\nu_i(x_i)~ d\nu_{u^c}(x_{u^c}) ~d\nu_{u\setminus i}(x_{u\setminus i})\\
&=& \displaystyle{\int \eta_v(x_v) \left(\int \eta_u(x_u) p_X(x) d\nu_i(x_i) ~d\nu_{u^c}(x_{u^c})\right)d\nu_{u\setminus i}(x_{u\setminus i})}\\
&=& 0 \qquad \textrm{ by assumption}
\end{array}
$$

\bigskip

Conversely, assume that ~$\exists~ u, ~\exists i \in u$ such that $\int \eta_u(x_u)p_X(x) d\nu_i (x_i)~ d\nu_{u^c}(x_{u^c}) \neq 0$, then,

\[
 \eta_v=\int \eta_u(x_u)p_X(x) d\nu_i (x_i)~ d\nu_{u^c}(x_{u^c}) \quad \textrm{with }v=u\setminus i
\]

and 

$$
\begin{array}{lll}
 \displaystyle{\int \eta_u(x_u) \eta_v(x_v)p_X(x)d\nu(x)} &=& \displaystyle{\int \eta_u(x_u) \left( \int \eta_u(x_u)p_X(x) d\nu_i ~ d\nu_{u^c}\right) p_X(x) d\nu_i ~d\nu_{u^c} ~d\nu_{u \setminus i}}\\
&=& \displaystyle{\int \!\!\! \left(\int \eta_u(x_u) p_X(x) d\nu_i(x_i) ~d\nu_{u^c}(x_{u^c})\right)^2 d\nu_{u \setminus i}(x_{u\setminus i})}\\
& > &0
\end{array}
$$

There is a contradiction, so that $\int \eta_u(x_u)p_X(x) d\nu_i (x_i)~ d\nu_{u^c}(x_{u^c})= 0 \quad \forall i \in u, \forall u$.\\
The second expression can be rewritten as :
\[
\int \eta_u(x_u)p_X(x) d\nu_i (x_i)~ d\nu_{u^c}(x_{u^c}) =\mathbb E(\eta_u/X_{u\setminus i}) \quad \forall i \in u, \forall u \in S
\]

Then, by Theorem \ref{theor1}, the minimization problem $(\mathcal P)$ admits a unique solution.

\end{preuve}


\subsection{Generalized sensitivity indices}\label{ap2}

\begin {preuve}{Proof of Proposition \ref{pro2}}

Under (\ref{c1}) and (\ref{c2}), Theorem \ref{theor1} states the existence and the uniqueness decomposition of $\eta$:

\[
 \eta(\mathbf X)=\sum_{u \in S}\eta_u(X_u), 
\]
 with $H_u^0 \perp H_v^0$, $\forall~v \subset u$. Therefore,

\[
 \mathbb E(\eta(X))= \mathbb E(\sum_{u \in S} \eta_u(X_u)) =\eta_0  
\]

and

$$
\begin{array}{lll}
 V(Y)=V(\eta(X))&=&\displaystyle{\mathbb E(\eta^2(X))-\eta_0^2} \\
&=& \displaystyle{\sum_{u \neq \emptyset} \mathbb E( \eta_u^2(X_u))+\sum_{u \neq v}\mathbb E(\eta_u(X_u)\eta_v(X_v))} \\
&=& \displaystyle{\sum_{u \neq \emptyset} V(\eta_u(X_u)) + \sum_{\substack{u \neq \emptyset \\ u\neq \mathcal P_p}}\sum_{\substack{v \neq \emptyset \\ u \nsubseteq v, v \nsubseteq u}} \mathbb E (\eta_u(X_u),\eta_v(X_v)) } \\
&=&\displaystyle{\sum_{u \neq \emptyset}\left[ V(\eta_u(X_u)) + \sum_{\substack{v\\ u \cap v \neq u,v}} \mbox{Cov}(\eta_u(X_u),\eta_v(X_v))\right]}
\end{array}
$$

Thus, (\ref{eq10}) holds, and equalities (\ref{eq11}) and (\ref{eq12}) follow obviously.

\end {preuve}


\section{Examples of distribution function}

\subsection{Boundedness of the inputs density function}

\label{annexp5}
\begin{preuve}{Proof of Proposition \ref{p5}}
Let $u\subset \mathcal P_p$, and $0<M_1\leq p_X\leq M_2$. As $ p_{X_u}$ and $ p_{X_u^c}$ are marginals, they are upper bounded by $M_2$.

As a consequence, $p_{X_u}p_{X_{u^c}}\leq\dfrac{M_2^2}{M_1}M_1\leq\dfrac{M_2^2}{M_1}p_X$, so that $ p_X\geq M_1M_2^{-2} p_{X_u}p_{X_{u^c}}$, with $0<M_1M_2^{-2}<1$.

\end{preuve}


\label{annexpro5}
\begin{preuve}{Proof of Example \ref{pro5}}

 \begin{itemize}
  \item $\nu$ is a product of measure as $\dfrac{d\nu}{d\nu_L}=\prod_{i=1}^p \nu_i(x_i)$, with $ \nu_i \sim N(m_i,\sigma_i^2)$. 
So 

\[
 \nu=\bigotimes_{i=1}^p \nu_i
\]

\item 
$p_X$ is given by

\begin{eqnarray}\label{densitygauss}
  p_X(\mathbf x)&=&\dfrac{dP_X}{d\nu}(\mathbf x)=\dfrac{dP_X}{d\nu_L}\times \dfrac{d\nu_L}{d\nu}(\mathbf x)\nonumber\\
 &=& \alpha + (1-\alpha)\left|\dfrac{\Sigma}{\Omega}\right|^{1/2}\exp{-\dfrac{1}{2}{}^t\!(\mathbf x-m)(\Omega^{-1}-\Sigma^{-1})(\mathbf x-m)}
\end{eqnarray}

 First, we have $ p_X(\mathbf x) \geq \alpha>0$. \\
 Further, the sufficient and necessary condition to have $p_X \leq M_2< \infty$ is to get $(\Omega^{-1}-\Sigma^{-1})$ positive definite.
 Indeed, if $(\Omega^{-1}-\Sigma^{-1})$ admits a negative eigenvalue, $p_X$ can not be bounded. 
Thus, $0<\alpha \leq p_X \leq M_2$ iff $(\Omega^{-1}-\Sigma^{-1})$ is positive definite.\\

 \end{itemize}
\end{preuve}


 \subsection{Examples of distribution of two inputs}

\begin{preuve}{Proof of Proposition \ref{copcond3}}

Condition (\ref{c4}) is immediate with Equation \ref{copula1}. Let prove that (\ref{c4}) is equivalent to (\ref{condcop}).\\

If (\ref{condcop}) holds, then $c(u,v) \geq M$. Conversely, we assume that $0<M < 1$, and
\[
\tilde C(u,v)=\frac{C(u,v)-Muv}{1-M}
\]
It is enough to show it is a copula : Obviously, $\tilde C(0,u)=C(u,0)=0$ and $\tilde C(1,u)=\tilde C(u,1)=u\quad \forall~u \in [0,1]$. 
By second order derivation, it comes that  $\tilde c(u,v)=\dfrac{c(u,v)-M}{1-M}$, so  $\tilde c(u,v)\geq 0$ by hypothesis (\ref{c4}).

\end{preuve}


\section{Estimation}\label{ap3}


\subsection{Model of $p=2$ input variables}\label{ap22}

\begin{preuve}{Proof of Proposition \ref{pro3}}

\begin{itemize}
 \item We first show first that ($\mathcal S$) admits an unique solution.\\

 Under (\ref{c1}) and (\ref{c2}), by Theorem \ref{theor1}, the decomposition of~ $\eta(\mathbf X)$ is unique and

\[
\eta(X_1,X_2)=\eta_0+\eta_1(X_1)+\eta_2(X_2)+\eta_{12}(X_1,X_2)
\]

with 

$$
\left\{
\begin{array}{l}
 \eta_0 \in H_{\emptyset} \\
\eta_i \in H_i^0 \perp H_{\emptyset}, \quad i=1,2 \\
\eta_{12} \in H_{12}^0 \perp H_i^0, \quad i=1,2, H_{12}^0 \perp H_{\emptyset}
\end{array}
\right.
$$

\bigskip

Thus,

 \begin{equation}\label{eq26}
\begin{pmatrix}
\id{Id} & 0 & 0 & 0\\
 0 & \id{Id} & P_{H_1^0} & 0\\
0 & P_{H_2^0} & \id{Id} & 0\\
0 & 0 &0 & \id{Id}
\end{pmatrix}
\begin{pmatrix}
\eta_0\\
 \eta_1\\
\eta_2\\
\eta_{12}
\end{pmatrix}
=
\begin{pmatrix}
P_{H_{\emptyset}}(\eta)\\
P_{H_1^0}(\eta)\\
P_{H_2^0}(\eta)\\
P_{H_{12}^0}(\eta)
\end{pmatrix}
\end{equation}

So $(\eta_0,\eta_1,\eta_2,\eta_{12})$ is solution of ($\mathcal S$). \\

Now, assume there exists an another solution of the system, say $(\widetilde\eta_0,\cdots,\widetilde{\eta}_{\mathcal P_p}) \in H_{\emptyset}\times\cdots\times H_{\mathcal P_p}^0$, then

$$
\begin{array}{lll}
 \left\{
\begin{array}{l}
\eta_0-\widetilde\eta_0 =0 \\
\eta_1-\widetilde\eta_1 + P_{H_1^0}(\eta_2-\widetilde\eta_2)=0 \\
P_{H_2^0}(\eta_1-\widetilde\eta_1)+\eta_2-\widetilde\eta_2=0\\
\eta_{12}-\widetilde\eta_{12}=0    
\end{array}
\right. & \Rightarrow & 
\left\{
\begin{array}{l}
\eta_0=\widetilde\eta_0 \\
 P_{H_1^0}(\eta_1-\widetilde\eta_1 +\eta_2-\widetilde\eta_2)=0 \\
P_{H_2^0}(\eta_1-\widetilde\eta_1+\eta_2-\widetilde\eta_2)=0\\
\eta_{12}=\widetilde\eta_{12}   
\end{array}
\right.\\
\\
&\Rightarrow & 
\left\{
\begin{array}{l}
\eta_0=\widetilde\eta_0 \\
\eta_1-\widetilde\eta_1 +\eta_2-\widetilde\eta_2 \in H_1^{0\perp} \cap H_2^{0\perp} \\
\eta_{12}=\widetilde\eta_{12}   
\end{array}
\right.
\end{array}
$$

As $\eta_1-\widetilde\eta_1 \in H_1^0$ and $\eta_2-\widetilde\eta_2 \in H_2^0$, it follows that

$$
\begin{array}{lll}
 \left\{
\begin{array}{l}
\langle\eta_1-\widetilde\eta_1 ,\eta_1-\widetilde\eta_1 +\eta_2-\widetilde\eta_2\rangle=0 \\
\langle\eta_2-\widetilde\eta_2 ,\eta_1-\widetilde\eta_1 +\eta_2-\widetilde\eta_2\rangle=0 \\
\end{array}
\right. & \Rightarrow & 
 \left\{
\begin{array}{l}
\|\eta_1-\widetilde\eta_1\|^2 +\langle\eta_1-\widetilde\eta_1,\eta_2-\widetilde\eta_2\rangle=0 \\
\|\eta_2-\widetilde\eta_2\|^2 +\langle\eta_1-\widetilde\eta_1 ,\eta_2-\widetilde\eta_2\rangle=0 \\
\end{array}
\right.\\
\\
&\Rightarrow & 
\|\eta_1-\widetilde\eta_1+\eta_2-\widetilde\eta_2\|^2=0 \\
&\Rightarrow& 
\eta_1-\widetilde\eta_1+\eta_2-\widetilde\eta_2=0
\end{array}
$$
 As $0$ can be uniquely decomposed into $H^0$ as a sum of zero, then,
\[
 \eta_1-\widetilde\eta_1=\eta_2-\widetilde\eta_2=0
\]

\item Let now compute
\begin{equation*} \label{eq14}
 \begin{pmatrix}
P_{H_{\emptyset}}(\eta)\\
P_{H_1^0}(\eta)\\
P_{H_2^0}(\eta)\\
P_{H_{12}^0}(\eta)
\end{pmatrix} 
\end{equation*}

First of all, it is obvious that the constant term $\eta_0=\mathbb E(\eta)$ and that $\eta_{12}$ is obtained by subtracting  $\eta$ with all other 
terms of the right of the decomposition.

Now, let us use the projector's property of embedded spaces. 
Indeed, as $H_i^0\subset H_i$, $\forall~i=1,2$, it comes

\[
 P_{H_i^0}(\eta)=P_{H_i^0}(P_{H_i}(\eta))=P_{H_i^0}[\underbrace{\mathbb{E}(\eta/X_i)}_{\varphi(X_i)}]
\]

$\varphi$ is a function of $X_i$, so it can be decomposed into the following expression :

\[
 \varphi(X_i)=\varphi_0+\varphi_i(X_i), \quad \varphi_0 \in H_{\emptyset}, ~\varphi_i \in H_i^0
\]

with $\varphi_0=\mathbb{E}(\varphi)=\mathbb E(\eta)$.\\
Hence, one can easily deduce $ P_{H_i^0}(\eta)$, $i=1,2$, as the term $\varphi_i=\mathbb{E}(\eta/X_i)-\mathbb{E}(\eta)$

We obtain

 \begin{equation}
 \begin{pmatrix}
P_{H_{\emptyset}}(\eta)\\
P_{H_1^0}(\eta)\\
P_{H_2^0}(\eta)\\
P_{H_{12}^0}(\eta)
\end{pmatrix} 
=
\begin{pmatrix}
\mathbb E(\eta)\\
\mathbb{E}(\eta/X_1)-\mathbb{E}(\eta)\\
\mathbb{E}(\eta/X_2)-\mathbb{E}(\eta)\\
\eta-\mathbb{E}(\eta/X_1)-\mathbb{E}(\eta/X_2)+\mathbb{E}(\eta)
\end{pmatrix}
\end{equation}

\end{itemize}

\end{preuve}

\subsection{Numerical procedure}\label{algo}

Gauss-Seidel algorithm requires the estimation of conditional expectation at each iteration. 
To do this, we use the local polynomial estimation with a leave-one-out technique. 
To considerably reduce time cost in application, the Sherman-Morrison formula is exploited.\\
We are going to review these methods for estimating $\mathbb E(Y/X=x)$, when $Y$ and $X$ are supposed to be real random variable.\\

The local polynomial estimation~\cite{fan} consists in approximating $m(x)=\mathbb E(Y/X=x)$ by a $q^{th}$-order polynomial fitted by a weight least squared estimation.

An explicit solution of $\hat m(x)$ is given by :

\begin{equation}\label{lpe}
 \hat m(x)={}^t(1~ 0\cdots 0)[{}^t\mathbb  X(x) D(x) \mathbb X(x))]^{-1}\cdot {}^t \mathbb X(x) D(x) \mathbb Y
\end{equation}

with

$$
\begin{array}{cc}
 \mathbb X(x) =
\begin{pmatrix}
 1 & X_1-x & \cdots & (X_1-x)^q \\
\vdots & \vdots & & \vdots \\
 1 & X_n-x & \cdots & (X_n-x)^q
\end{pmatrix}
& 
D(x)=
\begin{pmatrix}
K\left(\dfrac{X_1-x}{h}\right) & \cdots & 0 \\
& \ddots & \\
0 & \cdots & K\left(\dfrac{X_n-x}{h}\right)
\end{pmatrix}
\\
\mathbb Y=
\begin{pmatrix}
 Y_1\\
\vdots\\
Y_n
\end{pmatrix}

\end{array}
$$

 The leave-one-out technique on local estimation consists in estimating $m$ in every observation point $X_1,\cdots,X_n$, 
i.e. computing Equation (\ref{lpe}) when the $k^{th}$ line of matrices has been removed for estimating $\hat m(X_k)$.
It means that we would need to inverse ${}^t\mathbb  X_{-k}(x) D_{-k}(x) \mathbb X_{-k}(x)$ $n$ times, which is very expensive.
To avoid these expensive computations, Sherman and Morrison~\cite{SM} proposed a formula :

\begin{lem}
If $A$ is a square invertible matrix, and $u$, $v$ are vectors such that $1+{}^tvA^{-1}u \neq 0$, then

\begin{equation}
 (A+u{}^tv)^{-1}=A^{-1}-\dfrac{A^{-1}u{}^tvA^{-1}}{1+{}^tvA^{-1}u}
\end{equation}

\end{lem}

In our problem, set $S_n(x)={}^t\mathbb  X(x) D(x) \mathbb X(x)$. $S_n(x)$ can be rewritten as :

\begin{equation}
 S_n(x)=\sum_ {i=1}^n \Phi_i(x){}^t\Phi_i(x)=\underbrace{\sum_ {i\neq k} \Phi_i(x){}^t\Phi_i(x)}_{S_{-k}(x)}+\Phi_k(x){}^t\Phi_k(x)
\end{equation}

where
\[
 \Phi_i(x)={}^t(K(\dfrac{X_i-x}{h})\mathbb X_i(x)), \quad \mathbb X_i(x)=(1\cdots (X_i-x)^q),\quad \forall~i=1,\cdots,n
\]

Thus, $S_{-k}(x)$, corresponding to ${}^t\mathbb  X(x) D(x) \mathbb X(x)$ when the $k^{th}$ line has been removed, is of the form :

\[
 S_{-k}(x)=S_n(x)-\Phi_k(x){}^t\Phi_k(x), \quad \forall~k=1,\cdots,n
\]

The Sherman-Morrison formula gives :

\begin{equation}\label{invsm}
 S_ {-k}^{-1}(x)=S_n^{-1}(x)+\dfrac{S_n^{-1}(x)\Phi_k(x){}^t\Phi_k(x)S_n^{-1}(x)}{1-{}^t\Phi_k(x)S_n^{-1}(x)\Phi_k(x)}, \quad \forall~k=1,\cdots,n
\end{equation}

As $\hat m(X_k)={}^t(1~ 0\cdots 0) S_ {-k}^{-1}(X_k)\cdot {}^t \mathbb X_{-k}(X_k) D_{-k}(X_k) \mathbb Y$, $\forall~k$, it is faster to estimate $S_n^{-1}(x)$ and $\Phi_k(x)$ at each point of the design. \\
}

\section*{Acknowledgement}

This work has been supported by French National Research Agency (ANR) through COSINUS program (project COSTA-BRAVA
number ANR-09-COSI-015).

\bibliographystyle{imsart-number}
\bibliography{biblio}

\end{document}